\documentclass[]{amsart}
\usepackage{amssymb}

\usepackage{amsthm,amstext,amscd,latexsym,amsmath,amsfonts,amssymb}
\usepackage{amsmath}
\usepackage{graphicx}

\newtheorem{theorem}{Theorem}[section]

\newtheorem{definition}[theorem]{Definition}

\newcommand{\btimes}{\mathop{\LARGE\mbox{$\times$}}\nolimits}
\newcommand{\op}[1]{\mathop{\oplus}\limits_{\phantom{.}#1}}
\newcommand{\opp}[2]{\mathop{\oplus}\limits_{\phantom{.}#1}^{\phantom{.}#2}}

\newcommand{\pr}{\mbox{pr}_1\,}

\begin{document}
\title[Menger algebras of $n$-place functions]
      {Menger algebras of $n$-place functions}
\author{Wies{\l}aw A. Dudek}
\address{Institute of Mathematics and Computer Science\\
         Wroc{\l}aw University of Technology\\
         Wybrze\.ze Wyspia\'nskiego 27 \\
         50-370 Wroc{\l}aw, Poland}
\email{wieslaw.dudek@pwr.wroc.pl}
\author{Valentin S. Trokhimenko}
\address{Department of Mathematics\\
 Pedagogical University\\
 21100 Vinnitsa \\
 Ukraine}
\email{vtrokhim@gmail.com}

\begin{abstract}
It is a survey of the main results on abstract characterizations
of algebras of $n$-place functions obtained in the last $40$
years. A special attention is paid to those algebras of $n$-place
functions which are strongly connected with groups and semigroups,
and to algebras of functions closed with respect natural relations
defined on their domains.

 \end{abstract}
 \maketitle
 \centerline{\it{Dedicated to Professor K.P. Shum's 70th birthday}}
  \footnotetext{{\it 2010 Mathematics Subject Classification.}
           20N15}
\footnotetext{{\it Key words and phrases.} Menger algebra, algebra
of multiplace functions, representation, group-like Menger
algebra, diagonal semigroup}

\section{Introduction}

Every group is known to be isomorphic to some group of set
substitutions, and every semigroup is isomorphic to some semigroup
of set transformations. It accounts for the fact that the group
theory (consequently, the semigroup theory) can be considered as
an abstract study about the groups of substitutions (consequently,
the semigroups substitutions). Such approach to these theories is
explained, first of all, by their applications in geometry,
functional analysis, quantum mechanics, etc. Although the group
theory and the semigroup theory deal in particular only with the
functions of one argument, the wider, but not less important class
of functions remains without their attention -- the functions of
many arguments (in other words -- the class of multiplace
functions). Multiplace functions are known to have various
applications not only in mathematics itself (for example, in
mathematical analysis), but are also widely used in the theory of
many-valued logics, cybernetics and general systems theory.
Various natural operations are considered on the sets of
multiplace functions. The main operation is the superposition,
i.e., the operation which as a result of substitution of some
functions into other instead of their arguments gives a new
function. Although the algebraic theory of superpositions of
multiplace functions has not developed during a long period of
time, K.~Menger paid attention to abnormality of this state in the
middle of 40's in the previous century. He stated, in particular,
that the superposition of $n$-place functions, where $n$ is a
fixed natural number, has the property resembling the
associativity, which he called {\it super\-associativity }
\cite{134,138,139}. As it has been found later, this property
appeared to be fundamental in this sense that every set with
superassociative $(n+1)$-ary operation can be represented as a set
of $n$-place functions with the operation of superposition. This
fact was first proved by R.~M.~Dicker \cite{117} in 1963, and the
particular case of the given theorem was received by H.~Whitlock
\cite{157} in 1964.

The theory of algebras of multiplace functions, which now are
called {\it Menger algebras} (if the number of variables of
functions is fixed) or {\it Menger systems} (if the number of
variables is arbitrary), has been studied by (in alphabetic order)
M.~I.~Burtman \cite{6,8}, W.~A.~Dudek \cite{118} --
\cite{Dudtro4}, F.~A.~Gadzhiev \cite{18,18a,17}, L.~M.~Gluskin
\cite{23}, Ja.~Henno \cite{91} -- \cite{119b}, F.~A.~Ismailov
\cite{35b,35}, H.~L\"anger \cite{128,127,131a}, F.~Kh.~Muradov
\cite{44,44f}, B.~M.~Schein \cite{108,schtr}, V.~S.~Trokhimenko
(Trohimenko) \cite{61} -- \cite{76} and many others.

The first survey on algebras of multiplace functions was prepared
by B. M. Schein and V. S. Trokhimenko \cite{schtr} in 1979, the
first monograph (in Russian) by W. A. Dudek and V. S. Trokhimenko
\cite{Dudtro}. Extended English version of this monograph was
edited in 2010 \cite{Dudtro4}. This survey is a continuation of
the previous survey \cite{schtr} prepared in 1979 by B. M. Schein
and V. S. Trokhimenko.

\section{Basic definitions and notations}

An {\it $n$-ary relation} (or an {\it $n$-relation}) between
elements of the sets $A_1,A_2,\ldots,A_n$ is a subset  $\rho $ of
the Cartesian product  $A_1\times A_2\times\ldots\times A_n$. If
$A_1=A_2=\cdots=A_n$, then the $n$-relation $\rho$ is called {\it
homogeneous}. Later on we shall deal with $(n+1)$-ary relations,
i.e., the relations of the form $\rho\subset A_1\times
A_2\times\ldots\times A_n\times B$. For convenience we shall
consider such relations as binary relations of the form
$\rho\subset (A_1\times A_2\times\ldots\times A_n)\times B$ \ or
$\;\rho\subset\Big(\btimes\limits_{i=1}^{n}A_i\Big)\times B.$ In
the case of a homogeneous  $(n+1)$-ary relation we shall write
$\rho\subset A^n\times A$ or $\rho\subset A^{n+1}.$

Let $\,\rho\subset\Big(\btimes\limits_{i=1}^{n}A_i\Big)\times B\,$
be an $(n+1)$-ary relation, $\bar a= (a_1,\ldots ,a_n)$ an element
of $\,A_1\times A_2\times \ldots\times A_n$, $\,H_i\subset A_i$,
$i\in\{1,\ldots ,n\}=\overline{1,n}$, then \
$\rho\langle\bar{a}\rangle=\{b\in B\;|\;(\bar{a},b)\in\rho\}$ and
$$
\rho (H_1,\ldots ,H_n)=\bigcup\{\,\rho\langle\bar a\rangle\,
|\;\bar a\in H_1 \times H_2\times\ldots\times H_n\}.
$$
Moreover, let
\begin{eqnarray*}
&& \mbox{pr}_{1}\rho =\{\bar
a\in\btimes\limits_{i=1}^{n}A_i\;|\;(\exists b\in B)\,
(\bar a,b)\in\rho\},\nonumber\\[3pt]
&& \mbox{pr}_{2}\rho =\{b\in B\;|\;(\exists \bar a
\in\btimes\limits_{i=1}^{n}A_i)\, (\bar a,b)\in\rho\}.\nonumber
\end{eqnarray*}

To every sequence of $(n+1)$-relations $\,\sigma _1,\ldots ,\sigma
_n, \rho\,$ such that $\,\sigma _i\subset A_1\times\ldots \times
A_n\times B_i$, $\,i\in\overline{1,n}$,\, and $\,\rho\subset
B_1\times\ldots\times B_n \times C$, we assign an $(n+1)$-ary
relation
$$
\rho [\sigma_1\ldots \sigma _n]\subset A_1\times\ldots\times
A_n\times C,
$$
which is defined as follows:
$$
\rho [\sigma _1\ldots\sigma_n]=\{(\bar a,c)\;|\;(\exists \bar
b\,)\, (\bar a,b_1)\in\sigma_1\;\&\,\ldots\,\&\, (\bar
a,b_n)\in\sigma _n\;\&\,(\,\bar b,c)\in\rho\},
$$
where $\bar b =(b_1,\ldots ,b_n)\in B_1\times\ldots\times B_n$.
Obviously:
$$
\rho [\sigma _1\ldots\sigma _n](H_1,\ldots ,H_n)\subset \rho
(\sigma _1(H_1,\ldots ,H_n),\ldots,\sigma _n(H_1,\ldots ,H_n)\,),
$$
$$
\rho [\sigma _1\ldots\sigma _n][\chi _1\ldots \chi _n]\subset\rho
[\sigma _1[\chi _1\ldots \chi _n]\ldots \sigma _n[\chi _1\ldots
\chi _n]\,],
$$
where\footnote{It is clear that the symbol $\,\rho [\sigma
_1\ldots\sigma _n][\chi _1\ldots \chi _n]\,$ must be read as
$\,\mu [\chi _1\ldots \chi _n]$, where $\,\mu=\rho [\sigma
_1\ldots\sigma _n]$.} $\;\chi_i\subset A_1\times\ldots\times
A_n\times B_i$, \ $\sigma _i\subset B_1\times\ldots\times
B_n\times C_i$, \ $i=1,\ldots ,n$ and $\,\rho\subset
C_1\times\ldots\times C_n\times D$.

The $(n+1)$-operation $O:(\rho ,\sigma_1,\ldots ,\sigma _n)\mapsto
\rho [\sigma _1\ldots\sigma _n]$ defined as above is called a {\it
Menger superposition} or a {\it Menger composition} of relations.

Let
$$
\stackrel{n}{\triangle }_{A}=\{(\underbrace{a,\ldots ,a}_{n})\,|
\,a\in A\},
$$
then the homogeneous $(n+1)$-relation $\rho\subset A^{n+1}$ is
called
\begin{itemize}
\item  {\it reflexive } if \ $\stackrel{n+1}{\triangle }_{A}\subset\rho
$,
\item  {\it transitive } if \ $\rho [\rho\ldots \rho]\subset \rho
$,
\item  an $n$-{\it quasi-order} ({\it $n$-preorder}) if it is reflexive and
transitive. For $n=2$ it is a {\it quasi-order}.
\end{itemize}

An $(n+1)$-ary relation $\rho\subset A^n\times B$ is an {\it
$n$-place function} if it is one-valued, i.e.,
$$
(\forall \bar{a}\in A^n)(\forall b_1,b_2\in
B)\,(\,(\bar{a},b_1)\in\rho\,\;\&\,\;(\bar{a},b_2)\in\rho\,\longrightarrow\,
b_1=b_2).
$$

Any mapping of a subset of $A^n$ into $B$ is a {\it partial
$n$-place function}. The set of all such functions is denoted by
${\mathcal F}(A^n,B)$. The set of all {\it full $n$-place
functions} on $A$, i.e., mappings defined for every
$(a_1,\ldots,a_n)\in A^n$, is denoted by ${\mathcal T}(A^n,B)$.
Elements of the set ${\mathcal T}(A^n,A)$ are also called {\it
$n$-ary transformations of $A$}. Obviously ${\mathcal
T}(A^n,B)\subset{\mathcal F}(A^n,B)$. Many authors instead of a
full $n$-place function use the term an {\it $n$-ary operation}.

The superposition of $n$-place functions is defined by
\begin{equation}\label{1}
f[g_1\ldots g_n](a_1,\ldots,a_n) =
f(g_1(a_1,\ldots,a_n),\ldots,g_n(a_1,\ldots,a_n)),
\end{equation}
where $\,a_1,\ldots,a_n\in A$, $f,g_1,\ldots,g_n\in {\mathcal
F}(A^n,A)$. This superposition is an $(n+1)$--\\ ary operation
$\,\mathcal O\,$ on the set ${\mathcal F}(A^n,A)\,$ determined by
the formula $\,{\mathcal O}(f,g_1,\ldots,g_n) = f[g_1\ldots g_n]$.
Sets of $n$-place functions closed with respect to such
superposition are called {\it Menger algebras of $n$-place
functions} or {\it $n$-ary Menger algebras of functions}.

According to the general convention used in the theory of $n$-ary
systems, the sequence $\,a_i,a_{i+1},\ldots,a_j$, where
$i\leqslant j$, \ can be written in the abbreviated form as
$\,a_i^j$ \ (for \ $i>j$ \ it is the empty symbol). In this
convention (\ref{1}) can be written as
\[
f[g_1^n](a_1^n)=f(g_1(a_1^n),\ldots,g_n(a_1^n)).
\]
For $g_1=g_2=\ldots=g_n=g$ instead of $f[g_1^n]$ we will write
$f[g^n]$.

An $(n+1)$-ary groupoid $(G;o)$, i.e., a non-empty set $G$ with
one $(n+1)$-ary operation $o:G^{n+1}\rightarrow G$, is called a
{\it Menger algebra of rank $n$}, if it satisfies the following
identity (called the {\it superassociativity}):
\[
o(\,o(x,y_{1}^{n}),z_{1}^{n})=o(x,o(y_{1},z_{1}^{n}),o(y_{2},z_{1}^{n}),\ldots,
o(y_{n},z_{1}^{n})).
\]
A Menger algebra of rank $1$ is a semigroup.

Since a Menger algebra (as we see in the sequel) can be
interpreted as an algebra of $n$-place functions with a Menger
composition of such functions, we replace the symbol
$o(x,y_{1}^{n})$ by $x[y_{1}^{n}]$ or by $x[\bar{y}]$, i.e., these
two symbols will be interpreted as the result of the operation $o$
applied to the elements $x,y_{1},\ldots ,y_{n}\in G$.

In this convention the above superassociativity has the form
\[
x[y_{1}^{n}][z_{1}^{n}]=x[y_{1}[z_{1}^{n}]\ldots y_{n}[z_{1}^{n}]]
\]
or shortly
\[
x[\bar{y}][\bar{z}]=x[y_{1}[\bar{z}]\ldots y_{n}[\bar{z}]],
\]
where the left side can be read as in the case of functions, i.e.,
$x[y_{1}^{n}][z_{1}^{n}]=\big(x[y_1^n]\big)\,[z_1^n]$.

\medskip

\begin{theorem}\label{T21.2} {\bf (R. M. Dicker, \cite{117})}\newline
Every Menger algebra of rank $n$ is isomorphic to some Menger
algebra of full $n$-place functions.
\end{theorem}

Indeed, for every element $g$ of a Menger algebra $(G;o)$ of rank
$n$ we can put into accordance the full $n$-place function
$\lambda_{g}$ defined on the set $G\,^{\prime }=G\cup \{a,b\}$,
where $a$, $b$ are two different elements not belonging to $G$,
such that

\[ \lambda _{g}(x_{1}^{n})=
\left\{\begin{array}{cl}
g[x_{1}^{n}]&{\rm if } \ \ x_{1},\ldots ,x_{n}\in G, \\
g&{\rm if } \ \ x_{1}=\cdots =x_{n}=a, \\
b& {\rm in\; all\; other\; cases.}
\end{array}\right.
\]

Using the set $G^{\prime\prime}=G\cup\{a\}$, where $a\not\in G$,
and the partial $n$-place functions
\[ \lambda _{g}^{\prime }(x_{1}^{n})=
\left\{\begin{array}{cl}
g[x_{1}^{n}]&{\rm if } \ \ x_{1},\ldots ,x_{n}\in G, \\
g&{\rm if } \ \ x_{1}=\cdots =x_{n}=a,
\end{array}\right.
\]
we can see that every Menger algebra of rank $n$ is isomorphic to
a Menger algebra of partial $n$-place functions too \cite{108}.

\medskip

\begin{theorem} {\bf (J. Henno, \cite{97})}\newline
Every finite or countable Menger algebra of rank $n>1$ can be
isomorphically embedded into a Menger algebra of the same rank
generated by a single element.
\end{theorem}

\medskip

A Menger algebra $(G;o)$ containing {\it selectors}, i.e.,
elements $e_{1},\ldots ,e_{n}\in G$ such that
\[
x[e_{1}^{n}]=x \ \ \ \ {\rm and } \ \ \ \ e_{i}[x_{1}^{n}]=x_{i}
\]
for all $x,x_{i}\in G,$ $i=1,\ldots ,n$, is called {\it unitary}.

\begin{theorem} {\bf (V.S. Trokhimenko, \cite{61})}\newline
Every Menger algebra $(G;o)$ of rank $n$ can be isomorphically
embedded into a unitary Menger algebra $(G^{\ast};o^{\ast })$ of
the same rank with selectors $e_{1},\ldots ,e_{n}$ and a
generating set $G\cup\{e_{1},\ldots ,e_{n}\}$, where $e_{i}\not\in
G$ for all $i\in\overline{1,n}$.
\end{theorem}

\medskip

J. Hion \cite{100} and J. Henno \cite{98} have proven that Menger
algebras with selectors can be identified with some set of
multiplace endomorphisms of a universal algebra. Moreover E. Redi
proved in \cite{45} that a Menger algebra is isomorphic to the set
of all multiplace endomorphisms of some universal algebra if and
only if it contains all selectors.

W. N\"obauer and W. Philipp consider \cite{140} the set of all
one-place mappings of a fixed universal algebra into itself with
the Menger composition and proved  that for $n>1$ this algebra is
simple in the sense that it possesses no congruences other than
the equality and the universal relation \cite{142}.

\section{Semigroups}

The close connection between semigroups and Menger algebra was
stated already in 1966 by B.~M. Schein in his work \cite{108}. He
found other type of semigroups, which simply define Menger
algebras. Thus, there is the possibility to study semigroups of
such type instead of Menger algebras. But the study of these
semigroups is quite difficult, that is why in many questions it is
more advisable simply to study Menger algebras, than to substitute
them by the study of similar semigroups.

\begin{definition}\rm
Let $(G;o)$ be a Menger algebra of rank $n$. The set $G^{n}$
together with the binary operation $\ast$ defined by the formula:
\[
(x_{1},\ldots,x_{n})\ast(y_{1},\ldots,y_{n})=(x_{1}[y_{1}\ldots
y_{n}],\ldots,x_{n}[y_{1}\ldots y_{n}])
\]
is called the \textit{binary comitant} of a Menger algebra
$(G;o)$.
\end{definition}

An $(n+1)$-ary operation $o$ is superassociative if and only if
the operation $\ast$ is associative \cite{108}.

\begin{theorem}{\bf (B. M. Schein, \cite{schtr})}\newline
The binary comitant of a Menger algebra is a group if and only if
the algebra is of rank $1$ and is a group or the algebra is a
singleton.
\end{theorem}

It is evident that binary comitants of isomorphic Menger algebras
are isomorphic. However, as it was mentioned in \cite{108}, from
the isomorphism of binary comitants in the general case the
isomorphism of the corresponding Menger algebras does not follow.
This fact, due to necessity, leads to the consideration of binary
comitants with some additional properties such that the
isomorphism of these structures implies the isomorphism of initial
Menger algebras.

L.~M.~Gluskin observed (see \cite{23} and \cite{24}) that the sets
\[
M_{1}[G]=\{\bar{c}\in G^n
\,|\,x[\bar{y}][\bar{c}]=x[\bar{y}\ast\bar{c}\,]\;\;{\rm for\
all}\ x\in G,\;\bar{y}\in G^n \}
\]
and
\[
M_{2}[G]=\{\bar{c}\in G^n
\,|\,x[\bar{c}][\bar{y}]=x[\bar{c}\ast\bar{y}\,]\;\;{\rm for\
all}\ x\in G,\;\bar{y}\in G^n \}
\]
are either empty or subsemigroups of the binary comitant
$(G^{n};\ast)$. The set
\[
M_{3}[G]=\{a\in G
\,|\,a[\bar{x}][\bar{y}]=a[\bar{x}\ast\bar{y}\,]\;\;{\rm for\
all}\ \bar{x},\bar{y}\in G^n \}
\]
is either empty or a Menger subalgebra of $(G;o)$.

Let us define on the binary comitant $(G^{n};\ast)$ the
equivalence relations $\pi_{1},\ldots,\pi_{n}$ putting
\[
(x_{1},\ldots,x_{n})\equiv(y_{1},\ldots,y_{n})(\pi_{i})\longleftrightarrow
x_{i}=y_{i}
\]
for all $x_{i},y_{i}\in G,$ \ $i\in\overline{1,n}$. It is easy to
check, that these relations have the following properties:
\begin{itemize}
\item[$(a)$] for any elements $\bar{x}_{1},\ldots,\bar{x}_{n}\in G^{n}$ there is an
element $\bar{y}\in G^{n}$ such that $\bar{x}_{i}\equiv
\bar{y}(\pi_{i})$ for all $i\in\overline{1,n}$,
\item[$(b)$] if $\bar{x}\equiv\bar{y}(\pi_{i})$ for some $i\in\overline{1,n}$, then
$\bar{x}=\bar{y}$, where $\bar{x},\bar{y}\in G^{n}$,
\item[$(c)$] relations $\pi_{i}$ are {\it right regular}, i.e., $$
\bar{x}_{1}\equiv\bar{x}_{2}\longrightarrow\bar{x}_{1}\ast\bar{y}\equiv
\bar {x}_{2}\ast\bar{y}(\pi_{i})
$$
for all $\bar{x}_{1},\bar{x}_{2},\bar{y}\in G^{n}$, \
$i\in\overline{1,n}$,
\item[$(d)$] $\stackrel{n}{\triangle}_{G}$ is a {\it right ideal} of $(G^{n};\ast)$, i.e., $$
\bar{x}\in\,\stackrel{n}{\triangle}_{G}\,\wedge\;\bar{y}\in
G^{n}\longrightarrow\bar{x}\ast\bar{y}\in\,\stackrel{n}{\triangle}_{G}
$$
\item[$(e)$] for any $i\in\overline{1,n}$ every $\pi_{i}$-class contains precisely
one element from $\stackrel{n}{\triangle}_{G}$.
\end{itemize}

All systems of the type
$(G^{n};\ast,\pi_{1},\ldots,\pi_{n},\stackrel{n }{\triangle}_{G})$
will be called the \textit{rigged binary comitant of a Menger
algebra $(G;o)$}.
\begin{theorem}
Two Menger algebras of the same rank are isomorphic if and only if
their rigged binary comitants are isomorphic.
\end{theorem}

This fact is a consequence of the following more general theorem
proved in \cite{108}.

\begin{theorem}\label{T25.1}{\bf (B. M. Schein, \cite{108})}\newline
The system $(G;\cdot,\varepsilon_{1},\ldots,\varepsilon_{n},H)$,
where $(G;\cdot)$ is a semigroup,
$\varepsilon_{1},\ldots,\varepsilon_{n}$ are binary relations on
$G$ and $H\subset G$ is isomorphic to the rigged binary comitant
of some Menger algebra of rank $n$ if and only if
\begin{itemize}
\item[$1)$] for all $i=1,\ldots,n$ the relations $\varepsilon_{i}$ are
right regular equivalence relations and for any
$(g_{1},\ldots,g_{n})\in G^{n}$ there is exactly one $g\in G$ such
that $g_{i}\equiv g(\varepsilon_{i}) \ for\; all\;
i\in\overline{1,n},$
\item[$2)$] $H$ is a right ideal of a semigroup $(G;\cdot)$ and
for any $i\in\overline{1,n}$ every $\varepsilon_{i}$-class
contains exactly one element of $H$.
\end{itemize}
\end{theorem}

\medskip

In the literature the system of the type
$(G;\cdot,\varepsilon_{1},\ldots,\varepsilon _{n},H)$, satisfying
all the conditions of the above theorem, is called a
\textit{Menger semigroup of rank $n$}. Of course, a Menger
semigroup of rank 1 coincides with the semigroup $(G;\cdot)$.

So, the theory of Menger algebras can be completely restricted to
the theory of Menger semigroups. But we cannot use this fact. It
is possible that in some cases it would be advisable to consider
Menger algebras and in other -- Menger semigroups. Nevertheless,
in our opinion, the study of Menger semigroups is more complicated
than study of Menger algebras with one $(n+1)$-ary operation
because a Menger semigroup besides one binary operation contains
$n+1$ relations, which naturally leads to additional difficulties.

Let $(G;o)$ be a Menger algebra of rank $n$. Let us define on its
binary comitant $(G\,^{n};\ast)$ the unary operations
$\rho_{1},\ldots,\rho_{n}$ such that
\[
\rho_{i}(x_{1},\ldots,x_{n})=(x_{i},\ldots,x_{i})
\]
for any $x_{1},\ldots,x_{n}\in G$, $i\in\overline{1,n}$. Such
obtained system $(G^{n};\ast,\rho_{1},\ldots,\rho_{n})$ is called
the \textit{selective binary comitant} of Menger algebra $(G;o)$.

\begin{theorem}\label{T25.2}{\bf (B. M. Schein, \cite{108})}\newline
For the system $(G;\cdot,p_{1},\ldots,p_{n})$, where $(G;\cdot )$
is a semigroup and $p_{1},\ldots,p_{n}$ are unary operations on
it, a necessary and sufficient condition that this system be
isomorphic to the selective binary comitant of some Menger algebra
of rank $n$ is that the following conditions hold:
\begin{itemize}
\item[$1)$] $p_{i}(x)y=p_{i}(xy)$ for all $i\in\overline{1,n}$ and
$x,y\in G$,
\item[$2)$] $p_{i}\circ p_{j}=p_{j}$ for all $i,j\in\overline{1,n}$,
\item[$3)$] for every vector $(x_{1},\ldots,x_{n})\in G^{n}$ there is exactly
only one $g\in G$ such that $p_{i}(x_{i})=p_{i}(g)$ for all
$i\in\overline{1,n}$.
\end{itemize}
\end{theorem}

\medskip

The system $(G;\cdot,p_{1},\ldots,p_{n}) $ satisfying the
conditions of this theorem is called \textit{selective semigroups
of rank} $n$.

\begin{theorem} {\bf (B. M. Schein, \cite{108})}\newline
For every selective semigroup of rank $n$ there exists a Menger
algebra of the same rank with which the selective semigroup is
associated. This Menger algebra is unique up to isomorphism.
\end{theorem}

\medskip

These two theorems give the possibility to reduce the theory of
Menger algebras to the theory of selective semigroups. In this
way, we received three independent methods to the study of
superposition of multiplace functions: Menger algebras, Menger
semigroups and selective semigroups. A great number of papers
dedicated to the study of Menger algebras have been released
lately, but unfortunately the same cannot be said about Menger
semigroups and selective semigroups.

Defining on a Menger algebra $(G;o)$ the new binary operation
\[
x\cdot y=x[y\ldots y],
 \]
we obtain the so-called \textit{diagonal semigroup} $(G;\cdot )$.
An element $g\in G$ is called \textit{idempotent}, if it is
idempotent in the diagonal semigroup of $(G;o)$, i.e., if
$g[g^n]=g$. An element $e\in G$ is called a \textit{left
$($right$)$ diagonal unit} of a Menger algebra $(G;o)$, if it is a
left (right) unit of the diagonal semigroup of $(G;o)$, i.e., if
the identity $\,e[g^n]=g\,$ (respectively, $g[e^n]=g$) holds for
all $g\in G$. If $e$ is both a left and a right unit, then it is
called a {\it diagonal unit}. It is clear that a Menger algebra
has at most only one diagonal unit. Moreover, if a Menger algebra
has an element which is a left diagonal unit and an element which
is a right diagonal unit, then these elements are equal and no
other elements which are left or right diagonal units.

An $(n+1)$-ary groupoid $(G;o)$ with the operation
$o(x_{0}^{n})=x_{0}$ is a simple example of a Menger algebra of
rank $n$ in which all elements are idempotent and right diagonal
units. Of course, this algebra has no left units. In the Menger
algebra $(G;o_n)$, where $o_n(x_{0}^{n})=x_{n}$, all elements are
left diagonal units, but this algebra has no right diagonal units.

If a Menger algebra $(G;o)$ has a right diagonal unit $e$, then
its every element $c\in G$ satisfying the identity $e=e[c^n]$ is
also a right diagonal unit.

Non-isomorphic Menger algebras may have the same diagonal
semigroup. Examples are given in \cite{118}.

\begin{theorem}
A semigroup $(G;\cdot)$ is a diagonal semigroup of some Menger
algebra of rank $n$ only in the case when on $G$ can be defined an
idempotent $n$-ary operation $f$ satisfying the identity
$$f(g_{1},g_{2},\ldots,g_{n})\cdot g = f(g_{1}\cdot g,g_{2}\cdot
g,\ldots,g_{n}\cdot g).$$
\end{theorem}

\medskip

The operation of diagonal semigroup is in some sense distributive
with respect to the Menger composition. Namely, for all
$x,y,z_{1},\ldots,z_{n}\in G$ we have
\[(x\cdot y)[\bar{z}]= x\cdot y[\bar{z}]
\ \ \ \ {\rm and } \ \ \ \ x[\bar{z}]\cdot y = x[(z_{1}\cdot
y)\ldots (z_{n}\cdot y)].\]

In some cases Menger algebras can be completely described by their
diagonal semigroups. Such situation takes place, for example, in
the case of algebras of closure operations.

\begin{theorem}{\bf (V. S. Trokhimenko, \cite{Trokhim1})}\newline
A Menger algebra $(G;o)$ of rank $n$ is isomorphic to some algebra
of $n$-place closure operations on some ordered set if and only if
its diagonal semigroup $(G;\cdot)$ is a semilattice and
\[
x[y_{1}\ldots y_{n}]=x\cdot y_{1}\cdot\ldots\cdot y_{n}
\]
for any $x,y_{1},\ldots,y_{n}\in G$.
\end{theorem}

\medskip

A non-empty subset $H$ of a Menger algebra $(G;o)$ of rank $n$ is
called
\begin{itemize}
\item an {\it $s$-ideal\,} if \ $h[x_1^n]\in H$,

\item a {\it $v$-ideal\,} if \ $x[h_1^n]\in H$,

\item an {\it $l$-deal\,} if \ $x[x_1^{i-1},h_i,x_{i+1}^n]\in H$
\end{itemize}
for all $x,x_1,\ldots,x_n\in G$, $h,h_1,\ldots,h_n\in H$ and
$i\in\overline{1,n}$.

An $s$-ideal which is a $v$-ideal is called an {\it $sv$-ideal}. A
Menger algebra is {\it $s$-$v$-simple} if it possesses no proper
$s$-ideals and $v$-ideals, and {\it completely simple} if it
posses minimal $s$-ideals and $v$-ideals but has no proper
$sv$-ideals.

\begin{theorem} {\bf (Ja. N. Yaroker, \cite{114})}\newline
A Menger algebra is completely simple $(s$-$v$-simple$)$ if and
only if its diagonal semigroup is completely simple $($a group$)$.
\end{theorem}

\begin{theorem} {\bf (Ja. N. Yaroker, \cite{114})}\newline
A completely simple Menger algebra of rank $n$ can be decomposed
into disjoint union of $s$-$v$-simple Menger subalgebras with
isomorphic diagonal groups.
\end{theorem}

\medskip

Subalgebras obtained in this decomposition are classes modulo some
equivalence relation. This decomposition gives the possibility to
study connections between isomorphisms of completely simple Menger
algebras and isomorphisms of their dia\-gonal semigroups (for
details see \cite{114}).

If for $\bar{g}\in G^n$ there exists $x\in G$ such that
$g_{i}[x^n][\bar{g}]=g_{i}$ holds for each $i\in\overline{1,n}$,
then we say that $\bar{g}$ is a {\it $v$-regular} vector. The
diagonal semigroup of a Menger algebra in which each vector is
$v$-regular is a regular semigroup \cite{Trokhim2}.

An element $x\in G$ is an \textit{inverse element} for $\bar{g}\in
G^{n}$, if\ $x[\bar{g}][x^n]=x$\ and\ $g_{i}[x^n][\bar{g}]=g_{i}$
for all $i\in\overline{1,n}$. Every $v$-regular vector has an
inverse element \cite{Trokhim2}. Moreover, if each vector of a
Menger algebra $(G;o)$ is $v$-regular and any two elements of
$(G;o)$ are commutative in the diagonal semigroup of $(G;o)$, then
this semigroup is inverse.

\begin{theorem}\label{T22.5} {\bf (V. S. Trokhimenko,
\cite{Trokhim2})}\newline If in a Menger algebra $(G;o)$ each
vector is $v$-regular, then the diagonal semigroup of $(G;o)$ is a
group if and only if $(G;o)$ has only one idempotent.
\end{theorem}

\section{Group-like Menger algebras}
\setcounter{equation}{0}

A Menger algebra $(G;o)$ of rank $n$ in which the following two
equations
\begin{equation}
x[a_{1}\ldots a_{n}]=b,\label{sol-1}\\
\end{equation}
\begin{equation}
a_{0}[a_{1}\ldots a_{i-1}x_{i}a_{i+1}\ldots a_{n}]=b\label{sol-i}
\end{equation}
have unique solutions $x,x_{i}\in G$ for all $a_{0},a_{1}^n,b\in
G$ and some fixed $i\in\overline{1,n}$, is called {\it
$i$-solvable}. A Menger algebra which is $i$-solvable for every
$i\in\overline{1,n}$ is called {\it group-like}. Such algebra is
an $(n+1)$-ary quasigroup. It is associative (in the sense of
$n$-ary operations) only in some trivial cases (see \cite{Dud1,
Dud2, Dud3}).

A simple example of a group-like Menger algebra is the set of real
functions of the form
\[
f_{\alpha}(x_{1},\ldots,x_{n})=\frac{x_{1}+\cdots+x_{n}}{n}+\alpha.
\]

Every group-like Menger algebra of rank $n$ is isomorphic to some
group-like Menger algebra of {\it reversive} $n$-place functions,
i.e., a Menger algebra of functions $f\in\mathcal{F}(G^n,G)$ with
the property:
$$
f(x_1^{i-1},y,x_{i+1}^n)= f(x_1^{i-1},z,x_{i+1}^n)\,
\longrightarrow\, y=z
$$
for all $x_1^n,y,z\in G$ and $i=1,\ldots,n$.

The investigation of group-like Menger algebras was initiated by
H. Skala \cite{Sk}, the study of $i$-solvable Menger algebras by
W. A. Dudek (see \cite{Dud1}, \cite{Dud2} and \cite{118}). A
simple example of $i$-solvable Menger algebras is a Menger algebra
$(G;o_{i})$ with an $(n+1)$-ary operation
$o_{i}(x_{0},x_{1},\ldots,x_{n})=x_{0}+x_{i}$ defined on a
non-trivial commutative group $(G;+)$. It is obvious that the
diagonal semigroup of this Menger algebra coincides with the group
$(G;+)$. This algebra is $j$-solvable only for $j=i$, but the
algebra $(G;o)$, where
$o(x_{0},x_{1},\ldots,x_{n})=x_{0}+x_{1}+\ldots +x_{k+1}$ and
$(G;+)$ is a commutative group of the exponent $k\leqslant n-1$,
is $i$-solvable for every $i=1,\ldots,k+1$. Its diagonal semigroup
also coincides with $(G;+)$. Note that in the definition of
$i$-solvable Menger algebras one can postulate the existence of
solutions of \eqref{sol-1} and \eqref{sol-i} for all
$a_{0},\ldots, a_{n}\in G$ and some fixed $b\in G$ (see
\cite{118}). The uniqueness of solutions cannot be dropped, but it
can be omitted in the case of finite algebras.

\begin{theorem}{\bf (W. A. Dudek, \cite{118})}\label{D-T13}\newline
A finite Menger algebra $(G;o)$ of rank $n$ is $i$-solvable if and
only if it has a left diagonal unit $e$ and for all
$a_{1},\ldots,a_{n}\in G$ there exist $x,y\in G$ such that
 \[
x[a_{1}\ldots a_{n}]=a_{0}[a_{1}\ldots a_{i-1}y\,a_{i+1}\ldots
a_{n}]=e.
 \]
\end{theorem}

\medskip

A diagonal semigroup of an $i$-solvable Menger algebra is a group
\cite{Dud1}. The question when a given group is a diagonal group
of some Menger algebra is solved by the following two theorems
proved in \cite{118}.

\begin{theorem}\label{CD-P1}
A group $(G;\cdot)$ is a diagonal group of some $i$-solvable
Menger algebra if and only if on $G$ can be defined an $n$-ary
idempotent operation $f$ such that the equation $
f(a_{1}^{i-1},x,a_{i+1}^{n})=b $ has a unique solution for all
$a_{1}^{n},b\in G$ and the identity
$$
f(g_{1},g_{2},\ldots,g_{n})\cdot g = f(g_{1}\cdot g,g_{2}\cdot
g,\ldots,g_{n}\cdot g)
$$
is satisfied.
\end{theorem}

\medskip

\begin{theorem}\label{P23.3a}
A group $(G;\cdot)$ is the diagonal group of some $n$-solvable
Menger algebra of rank $n$ if and only if on $G$ can be defined an
$(n-1)$-ary operation $f$ such that $f(e,\ldots,e)=e$ for the unit
of the group $(G;\cdot)$ and the equation
\begin{equation}
f(a_{1}\cdot x,\ldots,a_{n-1}\cdot x)=a_{n}\cdot x\label{23.111}
\end{equation}
has a unique solution $x\in G$ for all $a_{1},\ldots,a_{n}\in G$.
\end{theorem}

\medskip
On the diagonal semigroup $(G;\cdot)$ of a Menger algebra $(G;o)$
of rank $n$ with the diagonal unit $e$ we can define a new
$(n-1)$-ary operation $f$ putting
\[f(a_{1},\ldots,a_{n-1})=e[a_{1}\ldots a_{n-1}e]
\]
for all $a_{1},...,a_{n-1}\in G$. The diagonal semigroup with such
defined operation $f$ is called a \textit{rigged diagonal
semigroup} of $(G;o)$. In the case when $(G;\cdot)$ is a group,
the operation $f$ satisfies the condition
\begin{equation}
a[a_{1}\ldots a_{n}]=a\cdot f(a_{1}\cdot
a_{n}^{-1},\ldots,a_{n-1}\cdot a_{n}^{-1})\cdot a_{n},
\label{23.11}
\end{equation}
where $a_{n}^{-1}$ is the inverse of $a_{n}$ in the group
$(G;\cdot)$.

\begin{theorem}{\bf (W. A. Dudek, \cite{118})}\newline
A Menger algebra $(G;o)$ of rank $n$ is $i$-solvable for some
$1\leqslant i<n$ if and only if in its rigged diagonal group
$(G;\cdot,f)$ the equation
\[
f(a_{1},\ldots,a_{i-1},x,a_{i+1},\ldots,a_{n-1})=a_{n}
\]
has a unique solution $x\in G$ for all $a_{1},\ldots,a_{n}\in G$.
\end{theorem}

\medskip

\begin{theorem}{\bf (W. A. Dudek, \cite{118})}\newline
A Menger algebra $(G;o)$ of rank $n$ is $n$-solvable if and only
if in its rigged diagonal group $(G;\cdot,f)$ for all
$a_{1},\ldots,a_{n-1}\in G$ there exists exactly one element $x\in
G$ satisfying \eqref{23.111}.
\end{theorem}

\medskip

\begin{theorem}{\bf (B. M. Schein,
\cite{schtr})}\label{T23.3}
\newline A Menger algebra $(G;o)$ of rank $n$ is group-like if and
only if on its diagonal semigroup $(G;\cdot)$ is a group with the
unit $e$, the operation $f(a_{1}^{n-1})=e[a_{1}^{n-1}e]$ is a
quasigroup operation and for all $a_{1},\ldots,a_{n-1}\in G$ there
exists exactly one $x\in G$ satisfying the equation
\eqref{23.111}.
\end{theorem}

\medskip

Non-isomorphic group-like Menger algebras may have the same
diagonal group \cite{Sk}, but group-like Menger algebras of the
same rank are isomorphic only in the case when their rigged
diagonal groups are isomorphic.

In \cite{Sk} conditions are considered under which a given group
is a diagonal group of a group-like Menger algebras of rank $n$.
This is always for an odd $n$ and a finite group. However, if both
$n$ and the order of a finite group are even, then a group-like
Menger algebra of rank $n$ whose diagonal group is isomorphic to
the given group need exists. If $n$ is even, then such algebra
exists only for finite orders not of the form $2p$, where $p$ is
an odd prime. There are no group-like Menger algebras of rank $2$
and finite order $2p$. The existence of group-like Menger algebras
of order $2p$ and even rank $n$ greater than $2$ is undecided as
yet.

\section{Representations by $n$-place functions}

Any homomorphism $P$ of a Menger algebra $(G;o)$ of rank $n$ into
a Menger algebra $({\mathcal{F}}(A^{n},A);O)$ of $n$-place
functions (respectively, into a Menger algebra
$(\mathfrak{R}(A^{n+1});O)$ of $(n+1)$-ary relations), where $A$
is an arbitrary set, is called a {\it representation of $(G;o)$ by
$n$-place functions} (respectively, {\it by $(n+1)$-ary
relations}). In the case when $P$ is an isomorphism we say that
this representation is {\it faithful}. If $P$ and $P_{i},\,i\in
I,$ are representations of $(G;o)$ by functions from
${\mathcal{F}}(A^{n},A)$ (relations from $\mathfrak{R}(A^{n+1})$)
and $P(g)=\bigcup_{i\in I}P_{i}(g)$ for every $g\in G$ $P$, then
we say that $P$ is the \textit{union} of the family $(P_{i})_{i\in
I}$. If $A=\bigcup_{i\in I}A_{i}$, where $A_{i}$ are pairwise
disjoint, then the union $\bigcup_{i\in I}P_{i}(g)$ is called the
\textit{sum} of $(P_{i})_{i\in I}$.

\begin{definition}\rm A {\it determining pair} of a
Menger algebra $(G;o)$ of rank $n$ is any pair $(\varepsilon,W)$,
where $\varepsilon$ is a partial equivalence on
$(G^{\ast};o^{\ast})$, $W$ is a subset of $G^{\ast}$ and the
following conditions hold:
\begin{itemize}
\item[$1)$] $G\cup\{e_{1},\ldots,e_{n}\}\subset{\rm pr}_{1}\varepsilon$,
where $e_{1},\ldots,e_{n}$ are the selectors of
$(G^{\ast};o^{\ast})$,
\item[$2)$] $e_{i}\not\in W$ for all $i=1,\ldots,n$,
\item[$3)$] $g[\varepsilon\langle e_{1}\rangle\ldots\varepsilon\langle
e_{n}\rangle]\subset\varepsilon\langle g\rangle$ for all $g\in G$,
\item[$4)$] $g[\varepsilon\langle g_{1}\rangle\ldots\varepsilon\langle
g_{n}\rangle]\subset\varepsilon\langle g[g_1\ldots g_n]\rangle$
for all $g,g_1,\ldots,g_n\in G$,
\item[$5)$] if $W\not =\varnothing$, then $W$ is an $\varepsilon$-class and
$W\cap G$ is an $l$-ideal of $(G;o)$.
\end{itemize}
\end{definition}

Let $(H_{a})_{a\in A}$ be the family of all $\varepsilon$-classes
indexed by elements of $A$ and distinct from $W$. Let $e_{i}\in
H_{b_{i}}$ for every $i=1,\ldots,n$,\ $A_{0}=\{a\in
A\,|\,H_{a}\cap G\not =\varnothing\}$,\
$\mathfrak{A}=A_{0}^{n}\cup\{(b_{1},\ldots,b_{n})\}$,\
$B=G^{n}\cup \{(e_{1},\ldots,e_{n})\}$. Every $g\in G$ is
associated with an $n$-place function $P_{(\varepsilon,W)}(g)$ on
$A$, which is defined by
\[
(\bar{a},b)\in P_{(\varepsilon,W)}(g)\longleftrightarrow\bar{a}\in
\mathfrak{A}\,\wedge\, g[H_{a_{1}}\ldots H_{a_{n}}]\subset H_{b},
\]
where $\bar{a}\in A^{n},$\ $b\in A$. The representation
$P_{(\varepsilon,W)}:g\to P_{(\varepsilon,W)}(g)$ is called the
{\it simplest representation} of $(G;o)$.

\begin{theorem}\label{T1.4.2}
Every representation of a Menger algebra of rank $n$ by $n$-place
functions is the union of some family of its simplest
representations.
\end{theorem}

The proof of this theorem can be found in \cite{Dudtro4} and
\cite{Dudtro4a}.

\medskip

With every representation $P$ of $(G;o)$ by $n$-place functions
($(n+1)$-ary relations) we associate the following binary
relations on $G$:

\[\begin{array}{lll}
& \zeta_{P}=\{(g_{1},g_{2})\,|\,P(g_{1})\subset P(g_{2})\}, \\[4pt]
& \chi_{P}=\{(g_{1},g_{2})\,|\,{\rm pr}_{1}P(g_{1})\subset{\rm pr}
_{1}P(g_{2})\}, \\[4pt]
& \pi_{P}=\{(g_{1},g_{2})\,|\,{\rm pr}_{1}P(g_{1})={\rm
pr}_{1}P(g_{2})\},  \\[4pt]
& \gamma_{P}=\{(g_{1},g_{2})\,|\,{\rm pr}_{1}P(g_{1})\cap{\rm pr}
_{1}P(g_{2})\not =\varnothing\},   \\[4pt]
& \kappa_{P}=\{(g_{1},g_{2})\,|\,P(g_{1})\cap P(g_{2})\not
=\varnothing \}, \\[4pt]
& \xi_{P}=\{(g_{1},g_{2})\,|\,P(g_{1})\circ\triangle_{{\rm {\small
pr}} _{1}P(g_{2})}=P(g_{2})\circ\triangle_{{\rm {\small
pr}}_{1}P(g_{1})}\}.
\end{array}
\]
A binary relation $\rho$ defined on $(G;o)$is {\it projection
representable} if there exists a representation $P$ of $(G;o)$ by
$n$-place functions for which $\rho=\rho_P$.

It is easy to see that if $P$ is a sum of the family of
representations $(P_{i})_{i\in I}$ then $\sigma_{P}=\bigcap_{i\in
I}\sigma_{P_{i}}$ for\label{SP} $\sigma\in\{\zeta,\chi ,\pi,\xi\}$
and $\sigma_{P}=\bigcup_{i\in I}\sigma_{P_{i}}$ for $\sigma
\in\{\kappa,\gamma\}$.

Algebraic systems of the form $(\Phi;O,\zeta_{\Phi},\chi_{\Phi})$,
$(\Phi;O,\zeta_{\Phi})$, $(\Phi;O,\chi_{\Phi})$ are called:
\textit{fundamentally ordered projection $($f.o.p.$)$ Menger
algebras, fundamentally ordered $($f.o.$)$ Menger algebras and
projection quasi-ordered $($p.q-o.$)$ Menger algebras}.

A binary relation $\rho$ defined on a Menger algebra $(G;o)$ of
rank $n\geqslant 2$ is
\begin{itemize}
\item {\it stable}, if \
$(x,y),(x_1,y_1),\ldots,(x_n,y_n)\in\rho\longrightarrow (x[x_1^n
],\,y[y_1^n]\,)\in\rho$,
\item {\it $l$-regular}, if \
$(x,y)\in\rho\longrightarrow (x[x_1^n ],y[x_1^n]\,)\in\rho$,
\item {\it $v$-negative}, if  \
$(x,\,u[x_1^{i-1},y,x_{i+1}^n]\,)\in\rho\longrightarrow(x,y)\in\rho$
\end{itemize}
\noindent for all $u,x,y,x_1,\ldots,x_n,y_1,\ldots,x_n\in G$ and
$i\in\overline{1,n}$.

\begin{theorem}\label{T2.1.1} {\bf (V. S. Trokhimenko, \cite{77})}\newline
An algebraic system $(G;o,\zeta,\chi)$, where $(G;o)$ is a Menger
algebra of rank $n$ and $\zeta$, $\chi$ are relations defined on
$G$, is isomorphic to a f.o.p. Menger algebra of $(n+1)$-relations
if and only if $\zeta$ is a stable order and $\chi$ is an
$l$-regular and $v$-negative quasi-order containing $\zeta$.
\end{theorem}

\medskip

From this theorem it follows that each stable ordered Menger
algebra of rank $n$ is isomorphic to some f.o. Menger algebra of
$\,(n+1)$--relations; each p.q-o. Menger algebra of
$(n+1)$-relations is a Menger algebra of rank $n$ with the
$l$-regular  and $v$-negative quasi-order.

Let $T_{n}(G)$ be the set of all polynomials defined on a Menger
algebra $(G;o)$ of rank $n$. For any binary relation $\rho$ on
$(G;o)$ we define the relation $\zeta (\rho)\subset G\times G$
putting $(g_{1},g_{2})\in\zeta (\rho)$ if and only if there exist
polynomials $t_{i}\in T_{n}(G)$, vectors $\bar{z}_{i}\in
B=G^{n}\cup\{\bar{e}\}$ and pairs
$(x_{i},y_{i})\in\rho\cup\{(e_{1},e_{1}),\ldots ,(e_{n},e_{n})\}$,
where $e_{1},\ldots ,e_{n}$ are the selectors from $(G^{\ast
};o^{\ast })$,  such that
\[
g_{1}=t_{1}(x_{1}[\bar{z}_{1}])\,\wedge\,
\bigwedge\limits_{i=1}^{m}\left(
t_{i}(y_{i}[\bar{z}_{i}])=t_{i+1}(x_{i+1}[\bar{z}_{i+1}])\right)
 \,\wedge \,t_{m+1}(y_{m+1}[\bar{z}_{m+1}])=g_{2}
\]
for some natural number $m$. Such defined relation $\zeta (\rho )$
is the least stable quasi-order on $(G;o)$ containing $\rho$.

\begin{theorem} \label{T2.1.3}{\bf (V. S. Trokhimenko, \cite{77})}\newline
An algebraic system $(G;o,\chi )$, where $(G;o)$ is a Menger
algebra of rank $n\geqslant 2$, $\chi \subset G\times G$, $\sigma
=\{(x,g[x^n])\,|\,x,g\in G\}$ is isomorphic to some p.q-o. Menger
algebra of reflexive $(n+1)$-relations if and only if $\zeta
(\sigma )$ is an antisymmetric relation, $\chi $ is an $l$-regular
$v$-negative quasi-order and $(t(x[\bar{y}]), t(g[x^n
][\bar{y}]))\in\chi$ for all \ $t\in T_{n}(G),$ \ $x,g\in G,$ \
$\bar{y}\in G^{n}\cup\{\bar{e}\}.$
\end{theorem}

\medskip

The system $(G;o,\zeta ,\chi )$, where $(G;o)$ is a Menger algebra
of rank $n$ and $\zeta ,\chi \subset G\times G$ satisfy all the
conditions of Theorem \ref{T2.1.1} can be isomorphically
represented by transitive $(n+1)$-relations if and only if for
each $g\in G$ we have $(g[g^n],g)\in\zeta$.

Analogously we can prove that each stable ordered Menger algebra
of rank $n$ satisfying the condition $(g[g^n], g)\in\zeta$, is
isomorphic to some f.o. Menger algebra of transitive
$(n+1)$-relations \cite{77}. Every idempotent stable ordered
Menger algebra of rank $n\geqslant 2$ in which
$(x,y[x^n])\in\zeta$ is isomorphic to some f.o. Menger algebra of
$n$-quasi-orders.

\begin{theorem}{\bf (V. S. Trokhimenko, \cite{61})}\newline
An algebraic system $(G;o,\zeta ,\chi)$, where $o$ is an
$(n+1)$-operation on $G$ and $\zeta ,\chi \subset G\times G$, is
isomorphic to a fundamentally ordered projection Menger algebra of
$n$-place functions if and only if $o$ is a superassociative
operation, $\zeta$ is a stable order, $\chi $ -- an $l$-regular
$v$-negative quasi-order containing $\zeta $, and for all
$i\in\overline{1,n}$ \ $u,g,g_{1},g_{2}\in G,$ \ $\bar{w}\in
G^{n}$ the following two implications hold
\begin{eqnarray*}
&&(g_{1},g), (g_{2},g)\in\zeta\ \& \
(g_{1},g_{2})\in\chi\longrightarrow (g_{1},g_{2})\in\zeta,\\[4pt]
&&(g_{1},g_{2})\in\zeta\ \& \ (g,g_{1}),(g,
u[\bar{w}|_{i}g_{2}])\in\chi\longrightarrow
(g,u[\bar{w}|_{i}g_{1}])\in\chi,
\end{eqnarray*}
where $u[\bar{w}|_ig]=u[w_1^{i-1},g,w_{i+1}^n]$.
\end{theorem}

\medskip

\begin{theorem}{\bf (V. S. Trokhimenko, \cite{61})}\newline
The necessary and sufficient condition for an algebraic system
$(G;o,\zeta )$ to be isomorphic to some f.o. Menger algebra of
$n$-place functions is that $o$ is a superassociative $(n+1)$-ary
operation and $\zeta $ is a stable order on $G$ such that
 \[
(x,y), (z,t_1(x)),(z,t_2(y))\in\zeta\longrightarrow
(z,t_2(x))\in\zeta
 \]
for all $x,y,z\in G$ and $t_1,t_2\in T_n(G)$.
\end{theorem}

\medskip

Replacing the last implication by the implication
 \[
(z,t_1(x), (z,t_1(y)),(z,t_2(y))\in\zeta\longrightarrow
(z,t_2(x))\in\zeta
 \]
we obtain a similar characterization of algebraic systems
$(G;o,\zeta)$ isomorphic to f.o. Menger algebras of reversive
$n$-place functions (for details see \cite{schtr} or
\cite{Dudtro4,Dudtro4a}).

\section{Menger algebras with additional operations}\setcounter{theorem}{0}
\setcounter{equation}{0}

Many authors investigate sets of multiplace functions closed with
respect to the Menger composition of functions and other naturally
defined operations such as, for example, the set-theoretic
intersection of functions considered as subsets of the
corresponding Cartesian product. Such obtained algebras are
denoted by $(\Phi;O,\cap )$ and are called \textit{$\cap $-Menger
algebras}. Their abstract analog is called a \textit{$
\curlywedge$-Menger algebra of rank} $n$ or a \textit{Menger
$\mathcal{P}$-algebra}.

An abstract characterization of $\cap $-Menger algebras is given
in \cite{73} (see also \cite{Dudtro6}).

\begin{theorem} \label{T3.1.1}{\bf (V. S. Trokhimenko, \cite{73})}\newline
For an algebra $(G;o,\curlywedge )$ of type $(n+1,2)$ the
following statements are true:

$(a)$ \ $(G;o,\curlywedge )$ is a $\curlywedge$-Menger algebra of
rank $n$, if and only if

$\rule{6.5mm}{0mm}(i)$ \ $o$ is a superassociative operation,

$\rule{6mm}{0mm}(ii)$ \ $(G;\curlywedge )$ is a semilattice, and
the following two conditions
\begin{eqnarray} &&(x\curlywedge y)[\bar{z}]=x[\bar{z}]\curlywedge
y[\bar{z}],  \label{3.4.1} \\[4pt]
&&t_{1}(x\curlywedge y\curlywedge z)\curlywedge
t_{2}(y)=t_{1}(x\curlywedge y)\curlywedge t_{2}(y\curlywedge
z)\label{3.4.2}
\end{eqnarray}

\hspace*{13mm}hold for all $\,x,y,z\in G,\ \bar{z}\in G^{n}$ and
$\,t_{1},t_{2}\in T_{n}(G)$,

$(b)$ \ if $\,n>1$, then $\,(G;o,\curlywedge )\,$ is isomorphic to
a $\,\cap$-Menger algebra of reversive

\hspace*{7mm}$n$-place functions if the operation $o$ is
superassociative, $\,(G;\curlywedge )\,$ is a semi--

\hspace*{7mm}lattice, $(\ref{3.4.1})$ and
\begin{eqnarray}
&&u[\bar{z}|_{i}(x\curlywedge y)]=u[\bar{z}|_{i}x]\curlywedge
u[\bar{z}|_{i}y],  \label{3.4.3} \\[4pt]
&&t_{1}(x\curlywedge y)\curlywedge t_{2}(y)=t_{1}(x\curlywedge
y)\curlywedge t_{2}(x)\label{3.4.4}
\end{eqnarray}

\hspace*{7mm}are valid for all $\,i\in\overline{1,n},\ u,x,y\in
G,\ \bar{z}\in G^{n},\ t_{1},t_{2}\in T_{n}(G)$.
\end{theorem}

\medskip

\begin{theorem} \label{T3.1.4}{\bf (B. M. Schein, V. S. Trokhimenko, \cite{76})}\newline
An algebra $(G;o,\curlyvee )$ of type $(n+1,2)$ is isomorphic to
some Menger algebra $(\Phi;O)$ of $n$-place functions closed under
the set-theoretic union of functions if and only if the operation
$o$ is superassociative, $(G;\curlyvee )$ is a semilattice with
the semillatice order $\leqslant$ such that $x[\bar{z}]\leqslant
z_{i}$ and
\begin{eqnarray}
&&(x\curlyvee y)[\bar{z}]=x[\bar{z}]\curlyvee y[\bar{z}],\label{3.4.20}\\
&&u[\bar{z}|_{i}(x\curlyvee y)]=u[\bar{z}|_{i}x]\curlyvee
u[\bar{z}|_{i}y],  \label{3.4.28} \\[4pt]
&&x\leqslant y\curlyvee u[\bar{z}|_{i}z]\longrightarrow x\leqslant
y\curlyvee u[\bar{z}|_{i}x],  \label{3.4.30}
\end{eqnarray}
for all $x,y,z\in G$, $\bar{z}\in G^{n},$ $u\in G\cup
\{e_1,\ldots,e_n\}$, $i\in\overline{1,n}$.
\end{theorem}

\medskip

\begin{theorem} \label{T3.1.5}{\bf (B. M. Schein, V. S. Trokhimenko, \cite{76})}\newline
An abstract algebra $(G;o,\curlywedge ,\curlyvee )$ of type
$(n+1,2,2)$ is isomorphic to some Menger algebra of $n$-place
functions closed with respect to the set-theoretic intersection
and union of functions if and only if the operation $o$ is
superassociative, $(G,\curlywedge ,\curlyvee )$ is a distributive
lattice, the identities $(\ref{3.4.1})$, $(\ref{3.4.20})$,
$(\ref{3.4.28})$ and
\begin{equation}
x[(y_{1}\curlywedge z_{1})\ldots (y_{n}\curlywedge
z_{n})]=x[\bar{y}]\curlywedge z_{1}\curlywedge \ldots \curlywedge
z_{n}\label{3.4.32}
\end{equation}
are satisfied for all $x,y\in G$, $\bar{y},\bar{z}\in G^{n},$
$u\in G\cup \{e_1,\ldots,e_n\}$, $i\in\overline{1,n}$.
\end{theorem}

\medskip

The restriction of an $n$-place function $f\in{\mathcal F}(A^n,B)$
to the subset $H\subset A^n$ can be defined as a composition of
$f$ and $\bigtriangleup_H=\{(\bar{a},\bar{a})|\bar{a}\in H\}$,
i.e., $f|_H=f\circ\bigtriangleup_H$. The {\it restrictive product}
$\rhd$ of two functions $f,g\in{\mathcal F}(A^n,B)$ is defined as
$$
f\rhd g=g\circ\bigtriangleup_{\pr}f.
$$

\begin{theorem} \label{T3.2.1}{\bf (V. S. Trokhimenko, \cite{65})}\newline
An algebra $(G;o,\blacktriangleright )$ of type $(n+1,2)$ is
isomorphic to a Menger algebra of $n$-place functions closed with
respect to the restrictive product of functions if and only if
$(G;o)$ is a Menger algebra of rank $n$, $(G;\blacktriangleright
)$ is an idempotent semigroup and the following three identities
hold:
\begin{eqnarray}
&&x[(y_{1}\blacktriangleright z_{1})\ldots
(y_{n}\blacktriangleright z_{n})]= y_{1}\blacktriangleright\ldots
\blacktriangleright y_{n}\blacktriangleright x[z_1\ldots z_n],
\rule{20mm}{0mm} \label{3.5.1} \\[4pt]
&&(x\blacktriangleright y)[\bar{z}]=x[\bar{z}]\blacktriangleright
y[\bar{z}],\label{3.5.2} \\[4pt]
&&\label{3.5.3}x\blacktriangleright y\blacktriangleright
z=y\blacktriangleright x\blacktriangleright z.
\end{eqnarray}
\end{theorem}

\medskip

An algebra $(G;o,\curlywedge ,\blacktriangleright )$ of type
$(n+1,2,2)$ is isomorphic to a Menger algebra of $n$-place
functions closed with respect to the set-theoretic intersection
and the restrictive product of functions, i.e., to
$(\Phi;O,\cap,\rhd)$, if and only if $(G;o,\blacktriangleright )$
satisfies the conditions of Theorem \ref{T3.2.1}, $(G;\curlywedge
)$ is a semilattice and the identities $(x\curlywedge
y)[\bar{z}]=x[\bar{z}]\curlywedge y[\bar{z}]$, $x\curlywedge
(y\blacktriangleright z)=y\blacktriangleright (x\curlywedge z)$,
$(x\curlywedge y)\blacktriangleright y=x\curlywedge y$ are
satisfied \cite{73}. Moreover, if $(G;o,\curlywedge
,\blacktriangleright )$ also satisfies $(\ref{3.4.3})$, then it is
isomorphic to an algebra $(\Phi;O,\cap,\rhd)$ of reversive
$n$-place functions.

More results on such algebras can be found in \cite{Dudtro4} and
\cite{Dudtro4a}.

\section{$(2,n)$-semigroups}\setcounter{theorem}{0}
\setcounter{equation}{0}

On $\mathcal{F}(A^{n},A)$ we can define $n$ binary compositions
$\op{1\,},\ldots ,\op{n}$ of two functions by putting
\[
(f\op{i\,}g)(a_{1}^n)=f(a_{1}^{i-1},g(a_{1}^{n}),a_{i+1}^{n}).
\]
for all $f,g\in\mathcal{F}(A^{n},A)$ and $a_1,\ldots,a_n\in A$.
Since all such defined compositions are associative, the algebra
$(\Phi; \op{1\,},\ldots ,\op{n})$, where
$\Phi\subset\mathcal{F}(A^{n},A)$, is called a
\textit{$(2,n)$-semigroup of $n$-place functions}.

The study of such compositions of functions for binary operations
was initiated by Mann \cite{Man} . Nowadays such compositions are
called {\it Mann's compositions} or {\it Mann's superpositions}.
Mann's compositions of $n$-ary operations, i.e., full $n$-place
functions, were described by T. Yakubov in \cite{Yak}. Abstract
algebras isomorphic to some sets of operations closed with respect
to these compositions are described in \cite{Sok}.

Menger algebras on $n$-place functions closed with respect to
Mann's superpositions are called {\it Menger $(2,n)$-semigroups}.
Their abstract characterizations are given in \cite{Dudtro1}.

Any non-empty set $G$ with $n$ binary operations defined on $G$ is
also called an abstract $(2,n)$-semigroup. For simplicity, for
these operations the same symbols as for Mann's compositions of
functions will be used. An abstract $(2,n)$-semigroup having a
representation by $n$-place function is called {\it
representable}.

Further, for simplicity, all expressions of the form $(\cdots
((x\op{i_{1}}y_{1})\op{i_{2}}y_{2})\cdots )\op{i_{k}}y_{k}$ are
denoted by $x\op{i_{1}}y_{1}\op{i_{2}}\cdots
\op{i_{k}}y_{k}$\label{sym96} or, in the abbreviated form, by
$x\op{i_{1}}^{i_{k}}y_{1}^{k}$. The symbol
$\mu_{i}(\op{i_{1}}^{i_{s}}x_{1}^{s})$ will be reserved for the
expression $x_{i_{k}}\!\op{i_{k+1}}^{i_{s}}\!x_{k+1}^{s}$, if
$i\neq i_{1}$, $\ldots$, $i\neq i_{k-1}$, $\,i=i_{k}$ for some
$k\in\{1,\ldots, s\}$. In any other case this symbol is empty. For
example, $\mu _{1}(\op{2}x\op{1\,}y\op{3}z)=y\op{3}z$, $\mu
_{2}(\op{2}x\op{1\,}y\op{3}z)=x\op{1\,}y\op{3}z$, $\mu
_{3}(\op{2}x\op{1\,}y\op{3}z)=z$. The symbol $\mu
_{4}(\op{2}x\op{1\,}y\op{3}z)$ is empty.

\begin{theorem}{\bf (W. A. Dudek, V. S. Trokhimenko, \cite{Dudtro2})}\newline
A $(2,n)$-semigroup $(G;\op{1\;},\ldots,\op{n})$ has a faithful
representation by partial $n$-place functions if and only if for
all $\,g,x_1,\ldots,x_s,y_1,\ldots,y_k\in G$ and
$\,i_1,\ldots,i_s,j_1,\ldots,j_k\in\overline{1,n}$ the following
implication
\begin{equation}\label{imp-14}
\bigwedge\limits_{i=1}^{n}\Big(\mu_i(\opp{i_1}{i_s}x_1^s)
=\mu_i(\opp{j_1}{j_k}y_1^k)\Big)\longrightarrow
g\opp{i_1}{i_s}x_1^s=g\opp{j_1}{j_k}y_1^k
\end{equation}
is satisfied.
\end{theorem}

\medskip

A $(2,n)$-semigroup $(G;\op{1\;},\ldots,\op{n})$ satisfying the
above implication has also a faithful representation by full
$n$-place functions.

Moreover, for any $(2,n)$-semigroup any its representation by
$n$-place functions is a union of some family of its simplest
representations \cite{Dudtro2}.

\begin{theorem}\label{T5}{\bf (W. A. Dudek, V. S. Trokhimenko, \cite{Dudtro2})}\newline
An algebraic system $(G;\op{1\,},\ldots,\op{n},\chi)$, where
$(G;\op{1\,},\ldots,\op{n})$ is a $(2,n)$-semigroup and $\chi$ is
a binary relation on $\,G$, is isomorphic to a p.q-o.
$(2,n)$-semigroup of partial $\,n$-place functions if and only if
the implication $\eqref{imp-14}$ is satisfied and $\chi$ is a
quasi-order such that
$(x\opp{i_{1}}{i_{s}}z_{1}^{s},\,\mu_{j}(\opp{i_{1}}{i_{s}}z_{1}^{s}))\in
\chi$ and $(x\op{i\,}z,\,y\op{i\,}z)\in\chi$ for all
$(x,y)\in\chi$.
\end{theorem}

\medskip

In the case of Menger $(2,n)$-semigroups the situation is more
complicated since conditions under which a Menger
$(2,n)$-semigroup is isomorphic to a Menger $(2,n)$-semi\-group of
$n$-place functions are not simple.

\begin{theorem}\label{T2}{\bf (W. A. Dudek, V. S. Trokhimenko, \cite{Dudtro1})}\newline
A Menger $(2,n)$-semigroup $(G;o,\op{1},\ldots,\op{n})$ is
isomorphic to a Menger $(2,n)$-semi\-group of partial $n$-place
functions if and only if it satisfies the implication
$\eqref{imp-14}$ and the following three identities
\begin{eqnarray*}\label{8}
&&x\op{i}y[z_1^n]=x[z_1^{i-1},y[z_1^n],z_{i+1}^n]\, ,\\[4pt]
 \label{9}
&&x[y_1^n]\op{i}z=x[y_1\op{i}z\ldots y_n\op{i}z]\, ,\\[4pt]
  \label{10}
&&x\opp{i_1}{i_s}y_1^s=x[\mu_1(\opp{i_1}{i_s}y_1^s)\ldots
\mu_n(\opp{i_1}{i_s}y_1^s)]\, ,
\end{eqnarray*}
where $\{i_1,\ldots,i_s\}=\{1,\ldots,n\}$.
\end{theorem}

\medskip

Using this theorem one can prove that any Menger $(2,n)$-semigroup
of $n$-place functions is isomorphic to some Menger
$(2,n)$-semigroup of full $n$-place functions (for details see
\cite{Dudtro1}). Moreover, for any Menger $(2,n)$-semigroup
satisfying all the assumptions of Theorem \ref{T2} one can find
the necessary and sufficient conditions under which the triplet
$(\chi,\gamma,\pi)$ of binary relations defined on this Menger
$(2,n)$-semigroup be a projection representable by the triplet of
relations $(\chi_P,\gamma_P,\pi_P)$ defined on the corresponding
Menger $(2,n)$-semigroup of $n$-place functions \cite{Dudtro3}.
Similar conditions were found for pairs $(\chi,\gamma)$,
$(\chi,\pi)$ and $(\gamma,\pi)$. But conditions under which the
triplet $(\chi,\gamma,\pi)$ will be a faithful projection
representable have not yet been found.

\setcounter{equation}{0}
\section{Functional Menger systems}\label{s36}

On the set ${\mathcal{F}}(A^{n},A)$ we also can consider $n$ unary
operations ${\mathcal{R}}_{1},\ldots ,{\mathcal{R}}_{n}$ such that
for every function $f\in {\mathcal{F}}(A^{n},A)$ \
${\mathcal{R}}_{i}f$ is the restriction of $n$-place projectors
defined on $A$ to the domain of $f$, i.e., $ {\mathcal{R}}_{i}f=f
\vartriangleright I_{i}^{n}\,$ for every $i=1,\ldots ,n,$ where
$I_{i}^{n}(a_1^n)=a_i\,$ is the $i$th $n$-place projector on $A$.
In other words, ${\mathcal{R}}_{i}f$ is such $n$-place function
from ${\mathcal{F}}(A^{n},A)$, which satisfies the conditions
\[{\rm pr}_{1}{\mathcal{R}}_{i}f ={\rm pr}_{1}f ,
\ \ \ \ \ \ \bar{a}\in {\rm pr}_{1}f \longrightarrow
{\mathcal{R}}_{i}f (\bar{a})=a_{i}
\]
for any $\bar{a}\in A^{n}$. Algebras of the form
$(\Phi;O,{\mathcal{R}}_{1},\ldots ,{\mathcal{R}}_{n})$, where
$\Phi\subset{\mathcal{F}}(A^{n},A)$, are called \textit{functional
Menger system of rank $n$.} Such algebras (with some additional
operations) were firstly studied in \cite{120} and \cite{151}. For
example, V. Kafka considered in \cite{120} the algebraic system of
the form $(\Phi;O,{\mathcal{R}}_{1},\ldots
,{\mathcal{R}}_{n},{\mathcal{L}},\subset )$, where
$\Phi\subset\mathcal{F}(A^{n},A)$ and ${\mathcal{L}}f
=\stackrel{n+1}\triangle_{\!{\rm pr}_{2}f}$ for every $f\in\Phi$.
Such algebraic system satisfies the conditions:
\begin{enumerate}
\item[$(A_1)$] \ $(\Phi;O)$ is a Menger algebra of functions,
\item[$(A_2)$] \ $\subset$ \ is an order on $\Phi$,
\item[$(A_3)$] \ the following identities are satisfied
\[\left\{\begin{array}{l}
f[\mathcal{R}_{1}f\ldots\mathcal{R}_{n}f]=f,\\[4pt]
\mathcal{R}_{i}(f[g_1\ldots g_n])=
\mathcal{R}_{i}((\mathcal{R}_{j}f)[g_1\ldots g_n]),\\[4pt]
(\mathcal{L}f)[f\ldots f]=f ,
\end{array}\right.
\]
\item[$(A_4)$] \ $\mathcal{L}(f[g_1\ldots g_n])\subset\mathcal{L}f\,$ for
all $\,f,g_1\ldots g_n\in\Phi$,
\item[$(A_5)$] \ $\Phi$ has elements $I_1,\ldots,I_n$ such that $\,f[I_1\ldots
I_n]=f\,$ and
\[\left\{\begin{array}{l}
g_i\subset I_i\longrightarrow\mathcal{R}_{i}g_i\subset g_i ,\\[4pt]
f\subset \bigcap\limits_{j=1}^{n}I_j\longrightarrow
\mathcal{L}f\subset f,\\[4pt]
g_k\subset I_k\longrightarrow g_k[f_1\ldots f_n]\subset f_k,\\[4pt]
\bigwedge\limits_{k=1}^{n}(g_k\subset I_k)\longrightarrow
f[g_1\ldots g_n]\subset f,\\[4pt]
\mathcal{R}_{i} g_j\subset\bigcap\limits_{k=1}^{n}\mathcal{R}_{i}
g_k\longrightarrow I_j[g_1\ldots g_n]=g_j ,\\[4pt]
f\subset h\longrightarrow f=h[p_1\ldots p_n] \ \ {\rm
for\;some}\;\;p_1\subset I_1,\ldots,p_n\subset I_n
\end{array}\right.
\]
for all $f,f_1,\ldots,f_n,g_1,\ldots,g_n\in\Phi$ and
$i,j,k\in\{1,\ldots,n\}$.
\end{enumerate}
It is proved in \cite{120} that for any algebraic system
$(G;o,R_{1},\ldots ,R_{n},L,\leqslant)$ satisfying the above
conditions there exists $\Phi\subset\mathcal{F}(A^{n},A)$ and an
isomorphism $(G;o)$ onto $(\Phi;O)$ which transforms the order
$\leqslant$ into the set-theoretical inclusion of functions.
However, $A_1 - A_5$ do not give a complete characterization of
systems of the form $(\Phi;O,{\mathcal{R}}_{1},\ldots
,{\mathcal{R}}_{n},{\mathcal{L}},\subset )$. Such characterization
is known only for systems of the form
$(\Phi;O,{\mathcal{R}}_{1},\ldots ,{\mathcal{R}}_{n})$.

\begin{theorem} \label{T3.3.2}{\bf (V. S. Trokhimenko, \cite{71})}\newline
A functional Menger system $(\Phi;O,{\mathcal{R}}_{1},\ldots
,{\mathcal{R}}_{n})$ of rank $n$ is isomorphic to an algebra
$(G;o,R_{1},\ldots ,R_{n})$ of type $(n+1,1,\ldots ,1)$, where
$(G;o)$ is a Menger algebra, if and only if for all
$i,k\in\overline{1,n}$ it satisfies the identities
\begin{eqnarray}
&&\hspace*{-15mm}x[R_{1}x\ldots R_{n}x]=x,\label{3.6.5}  \\[4pt]
&&\hspace*{-15mm}x[\bar{u}|_{i}z][R_{1}y\ldots R_{n}y]=
x[\bar{u}|_{i}\,z[R_{1}y\ldots R_{n}y]],\label{3.6.6}  \\[4pt]
&&\hspace*{-15mm}R_{i}(x[R_{1}y\ldots R_{n}y])=
(R_{i}x)[R_{1}y\ldots R_{n}y],\label{3.6.7} \\[4pt]
&&\hspace*{-15mm}x[R_{1}y\ldots R_{n}y][R_{1}z\ldots R_{n}z]=
x[R_{1}z\ldots R_{n}z][R_{1}y\ldots R_{n}y],\label{3.6.8} \\[4pt]
&&\hspace*{-15mm}R_{i}(x[\bar{y}])=
R_{i}((R_{k}x)[\bar{y}]), \label{3.6.9}\\[4pt]
&&\hspace*{-15mm}(R_{i}x)[\bar{y}]=y_{i}[R_{1}(x[\bar{y}])\ldots
R_{n}(x[\bar{y}])].\label{3.6.10}
\end{eqnarray}
\end{theorem}

\medskip

It is interesting to note that defining on $(G;o,R_{1},\ldots
,R_{n})$ a new operation $\blacktriangleright$ by putting
$x\blacktriangleright y=y[R_{1}x\ldots R_{n}x]$ we obtain an
algebra $(G;o,\blacktriangleright )$ isomorphic to a Menger
algebra of $n$-place functions closed with respect to the
restrictive product of functions \cite{Dudtro4a}.

\medskip

Another interesting fact is
\begin{theorem} \label{T3.3.3}{\bf (W. A. Dudek, V. S. Trokhimenko, \cite{Dudtro})}\newline
An algebra $(G;o,\curlywedge,R_{1},\ldots,R_{n})$ of type
$(n+1,2,1,\ldots ,1)$ is isomorphic to some functional Menger
$\cap$-algebra of $n$-place functions if and only if
$(G;o,R_{1},\ldots,R_{n})$ is a functional Menger system of rank
$n$, $(G;\curlywedge )$ is a semilattice, and the identities
\begin{eqnarray*}
&&(x\curlywedge y)[z_1\ldots z_n]= x[z_1\ldots z_n]\curlywedge
y[z_1\ldots z_n],\\[4pt]
&&(x\curlywedge y)[R_1z\ldots R_nz]=x\curlywedge y[R_1z\ldots
R_nz],\\[4pt]
&&x[R_1(x\curlywedge y)\ldots R_n(x\curlywedge y)]=x\curlywedge y
 \end{eqnarray*}
are satisfied.
\end{theorem}

\medskip

Now we present abstract characterizations of two important sets
used in the theory of functions. We start with the set containing
functions with the same fixed point.

\begin{definition} A non-empty subset $H$ of $G$ is called a \textit{stabilizer} of a
functional Menger system $(G;o,R_1,\ldots,R_n)$ of rank $n$ if
there exists a representation $P$ of $(G;o,R_1,\ldots,R_n)$ by
$n$-place functions $f\in{\mathcal F}(A^n,A)$ such that
$$
H=H^a_P=\{g\in G\,|\,P(g)(a,\ldots,a)=a\}
$$
for some point $a\in A$ common for all $g\in H$.
\end{definition}

\begin{theorem}{\bf (W. A. Dudek, V. S. Trokhimenko, \cite{Dudtro5})}\newline
A non-empty $H\subset G$ is a stabilizer of a functional Menger
system $(G;o,R_1,\ldots,R_n)$ of rank $n$ if and only if there
exists a subset $U$ of $\,G$ such that
\begin{eqnarray}
&&\hspace*{-15mm} \label{f-0}H\subset U, \ \ R_iU\subset H, \ \
R_i(G\!\setminus\! U)\subset G\!\setminus\!U,\\
&&\hspace*{-15mm}\label{f-1}
x,y\in H\;\,\&\;\,t(x)\in U\longrightarrow t(y)\in U, \\
&&\hspace*{-15mm} \label{f-2} x=y[R_1x\ldots R_nx]\in U\;\,\&\;\,
u[\bar{w}\,|_iy]\in H\longrightarrow u[\bar{w}\,|_ix]\in H, \\
&&\hspace*{-15mm} \label{f-3} x=y[R_1x\ldots R_nx]\in
U\;\,\&\;\,u[\bar{w}\,|_iy]\in U \longrightarrow
u[\bar{w}\,|_ix]\in U,\\
&&\hspace*{-15mm}y\in H\;\,\&\;\,x[y\ldots y]\in H\longrightarrow x\in H,\\
&&\hspace*{-15mm}x,y\in H\;\,\&\;\,t(x)\in H\longrightarrow t(y),\\
&&\hspace*{-15mm}x\in H\longrightarrow x[x\ldots x]\in
H\label{f-13}
\end{eqnarray}
for all $x,y\in G,\,$ $\bar{w}\in G^{\,n},\,$ $t\in T_n(G)$ and
$i\in\overline{1,n}$, where the symbol $u[\bar{w}\,|_i\ ]$ may be
empty.\footnote{ \ If $u[\bar{w}\,|_i\ ]$ is the empty symbol,
then $u[\bar{w}\,|_ix]$ is equal to $x$.}
\end{theorem}

\medskip

In the case of functional Menger algebras isomorphic to functional
Menger $\cap$-algebras of $n$-place functions stabilizers have
simplest characterization. Namely, as it is proved in
\cite{Dudtro5}, the following theorem is valid.

\begin{theorem}
A non-empty subset $H$ of $\,G$ is a stabilizer of a functional
Menger algebra $(G;o,\curlywedge,R_1,\ldots,R_n)$ of rank $n$ if
and only if $H$ is a subalgebra of $(G;o,\curlywedge)$,
$\,R_iH\subset H$ for every $i\in\overline{1,n}$, and
$$
x[y_1\ldots y_n]\in H\longrightarrow x\in H
$$
for all $x\in G,$ $y_1,\ldots,y_n\in H$.
\end{theorem}

\medskip

Another important question is of abstract characterizations of
stationary subsets of Menger algebras of $n$-place functions. Such
characterizations are known only for some types of such algebras.

As it is well known, by the \textit{stationary subset} of
$\Phi\subset\mathcal{F}(A^n,A)$ we mean the set
 \[
\mathbf{St}(\Phi)= \{f\in\Phi\,|\,(\exists a\in
A)\,f(a,\ldots,a)=a\}.
 \]

Note that in a functional Menger $\cap$-algebra
$(\Phi;\mathcal{O},\cap,\mathcal{R}_1,\ldots,\mathcal{R}_n)$ of
$n$-place functions without a zero $\mathbf{0}$ the subset
$\mathbf{St}(\Phi)$ coincides with $\Phi$. Indeed, in this algebra
$f\neq\emptyset$ for any $f\in\Phi$. Consequently, $f\cap
g[f\ldots f]\neq\emptyset$ for all $f,g\in\Phi$. Hence, $g[f\ldots
f](\bar{a})=f(\bar{a})$, i.e.,
$g(f(\bar{a}),\ldots,f(\bar{a}))=f(\bar{a})$ for some $\bar{a}\in
A^n$. This means that $g\in\mathbf{St}(\Phi)$. Thus,
$\Phi\subseteq\mathbf{St}(\Phi)\subseteq\Phi$, i.e.,
$\mathbf{St}(\Phi)=\Phi$. Therefore, below we will consider only
functional Menger $\cap$-algebras with a zero.

\begin{definition} \label{D-2}\rm
A non-empty subset $H$ of $G$ is called a \textit{stationary
subset} of a functional Menger algebra
$(G;o,\curlywedge,R_{1},\ldots,R_{n})$ of rank $n$ if there exists
its faithful representation $P$ by $n$-place functions such that
\begin{equation*}
 g\in H\longleftrightarrow P(g)\in\mathbf{St}(P(G))
\end{equation*}
for every $g\in G$, where $P(G)=\{P(g)\,|\,g\in G\}$.
\end{definition}

\begin{theorem}{\bf (W. A. Dudek, V. S. Trokhimenko, \cite{Dudtro6})}\newline
A non-empty subset $H$ of $G$ is a stationary subset of a
functional Menger algebra $(G;o,\curlywedge,R_{1},\ldots,R_{n})$
with a zero $\mathbf{0}$ if and only if
\begin{eqnarray*}
&&\hspace*{-60mm}x[x\ldots x]\curlywedge x\in H,\\[4pt]
&&\hspace*{-60mm}  z[y\ldots y]=y\in H\longrightarrow z\in H, \\[4pt]
&&\hspace*{-60mm}  z[y\ldots y]\curlywedge
y\neq\mathbf{0}\longrightarrow z\in H,\\[4pt]
&&\hspace*{-60mm} \mathbf{0}\not\in H\longrightarrow
R_i\mathbf{0}=\mathbf{0},
\end{eqnarray*}
for all $x\in H,$ $\,y,z\in G$ and $i\in\overline{1,n}$.
\end{theorem}

\medskip

Conditions formulated in the above theorem are not identical with
the conditions used for the characterization of stationary subsets
of restrictive Menger $\mathcal P$-algebras (see Theorem 8 in
\cite{Trokh3}). For example, the implication
$$
y[R_1x\ldots R_nx]=x\in H\longrightarrow y\in H
$$
is omitted. Nevertheless, stationary subsets of functional Menger
$\curlywedge$-algebras with a zero have the same properties as
stationary subsets of restrictive Menger $\mathcal P$-algebras.

More results on stationary subsets and stabilizers in various
types of Menger algebras can be found in \cite{Dudtro4} and
\cite{Dudtro4a}. In \cite{Dudtro4a} one can find abstract
characterizations of algebras of vector-valued functions,
positional algebras of functions, Mal'cev--Post iterative algebras
and algebras of multiplace functions partially ordered  by various
types of natural relations connected with domains of functions.

\end{document}